\documentstyle[12pt]{article}

\newtheorem{theorem}{Theorem}
\newtheorem{lemma}{Lemma}
\newtheorem{definition}{Definition}

\newcommand{\be}{\begin{equation}}
\newcommand{\ee}{\end{equation}}

\begin{document}
\bibliographystyle{plain}

\thispagestyle{empty}
\setcounter{page}{0}

\vspace {2cm}

{\Large Guszt\'av Morvai and   Benjamin Weiss: }

\vspace {1cm}

{\LARGE Forecasting for stationary binary time series.}

\vspace {1cm}

{\Large  Acta Appl. Math.  79  (2003),  no. 1-2, 25--34.}

\vspace {2cm}

\begin{abstract}
The forecasting problem for a  stationary and ergodic binary time series $\{X_n\}_{n=0}^{\infty}$ is to 
estimate the probability that $X_{n+1}=1$ based on the observations 
$X_i$, $0\le i\le n$ without prior knowledge of the distribution of the process $\{X_n\}$. 
It is known that this is not possible if one estimates at all values of $n$. 
 We present a simple procedure which will attempt to make such a prediction 
 infinitely often at carefully selected stopping times chosen 
by the algorithm. We show that the proposed procedure is consistent under certain 
conditions, and we estimate the growth rate of the stopping times. 
\end{abstract}
%\keywords{Nonparametric estimation, stationary processes}
%\classification{Mathematics Subject Classifications (2000)}{62G05, 60G25, 60G10}
%\end{opening}

\pagebreak

\section{Introduction}

T. Cover \cite{Cover75} posed two fundamental problems concerning 
estimation for stationary and ergodic binary  time series $\{ X_n\}_{n=-\infty}^{\infty}$.
(Note that a stationary time series $\{ X_n\}_{n=0}^{\infty}$ 
can be extended to be a two sided stationary time series 
$\{ X_n\}_{n=-\infty}^{\infty}$.) 
Cover's first problem was on  backward estimation. 

\smallskip
\noindent
{\bf Problem 1} {\it Is there an estimation scheme $f_{n+1}$ for the value 
$P(X_1=1|X_{-n},\dots,X_{0})$ 
such that $f_{n+1}$ depends solely on the observed data 
segment $(X_{-n},\dots,X_{0})$ 
and 
$$
\lim_{n\to\infty} |f_{n+1}(X_{-n},\dots,X_{0})-P(X_{1}=1|X_{-n},\dots,X_{0})|=0
$$
almost surely for all stationary and ergodic binary   
time series $\{ X_n\}_{n=-\infty}^{\infty}$? 
}

\smallskip
\noindent
This problem was solved by Ornstein \cite{Ornstein78} by constructing such a scheme. (See also Bailey \cite{Bailey76}.)
Ornstein's scheme is not 
a simple one and  the proof of consistency is rather sophisticated. 
A much simpler scheme and proof of consistency were provided by Morvai, Yakowitz, Gy\"orfi \cite{MoYaGy96}. (See also Weiss \cite{Weiss00}.)

\smallskip
\noindent
Cover's second problem was on forward estimation (forecasting).  

\smallskip
\noindent
{\bf Problem 2} {\it Is there an estimation scheme $f_{n+1}$ for the value 
$P(X_{n+1}=1|X_0,\dots,X_n)$ such that 
$f_{n+1}$ depends solely on the data segment $(X_0,\dots, X_n)$ and 
$$
\lim_{n\to\infty} |f_{n+1}(X_0,\dots,X_n)-P(X_{n+1}=1|X_0,\dots,X_n)|=0
$$
almost surely for all stationary and ergodic binary    
time series $\{ X_n\}_{n=-\infty}^{\infty}$? }

\smallskip
\noindent
This problem was answered by Bailey \cite{Bailey76} in a negative way, that is, he showed that there is no such scheme. 
(Also see Ryabko \cite{Ryabko88}, Gy\"orfi,  Morvai, Yakowitz \cite{GYMY98}
 and  Weiss \cite{Weiss00}.) 
Bailey used the technique of cutting and stacking developed by Ornstein \cite{Ornstein74} 
(see also Shields \cite{Shields91}). Ryabko's construction was  based on a function of 
an infinite state Markov-chain. 
This negative result can be interpreted as follows. Consider a market analyst  whose task
 it is to predict the probability of the event 'the price of a certain share will go up tomorrow' given the 
observations up to the present day. Bailey's result says that the difference between the
 estimate and the true conditional probability cannot   eventually be small for all 
stationary and ergodic market processes. The difference will be big infinitely often. 
These results show that there is a great  difference between Problems~1 and~2.  Problem~1 was 
addressed by Morvai, Yakowitz, Algoet \cite{MoYaAl97} and a very simple estimation scheme was 
given which satisfies the statement in Problem~1 {\it in probability} instead of 
{\it almost surely}. 
However, for the class of all stationary and ergodic binary Markov-chains of some finite  order Problem~2 can be solved.
 Indeed,  
if the time series is a Markov-chain of some finite (but unknown) order, we can estimate the order 
(e.g. as in Csisz\'ar, Shields 
\cite{CsSh00}) and  count frequencies of  blocks with length equal to the order.

\smallskip
\noindent Let ${\cal X}^{*-}$ be the set of all one-sided  binary sequences, that is, 
$${\cal X}^{*-} =\{ (\dots,x_{-1},x_0): x_i\in \{0,1\} \ \ \mbox{for all $-\infty<i\le 0$}\}.$$

\noindent
Let $d(\cdot,\cdot)$ be the Hamming distance (that is for $x,y\in \{0,1\}$, 
$d(x,y)=0$ 
if and only if $x=y$ and $d(x,y)=1$ otherwise), and define the distance on sequences
 $(\dots,x_{-1},x_{0},)$ and $(\dots,y_{-1},y_{0})$  as follows. Let 
\begin{equation}
\label{defdistance}
d^*((\dots,x_{-1},x_{0}),(\dots,y_{-1},y_{0}))=
\sum_{i=0}^{\infty} 2^{-i-1} d(x_{-i},y_{-i}).
\end{equation} 
(For details see Gray~\cite{Gr88} p. 51. )

\smallskip
\noindent
\begin{definition}
The conditional probability $P(X_1=1|\dots,X_{-1},X_{0})$
is almost surely continuous if for some set $C\subseteq {\cal X}^{*-}$  
which has probability one   
the conditional probability $P(X_1=1|\dots,X_{-1},X_{0})$ restricted to this set $C$ 
is continuous with respect to metric $ d^*(\cdot,\cdot)$ in (\ref{defdistance}). 
\end{definition}

\noindent
We note that  from the proof of Ryabko~\cite{Ryabko88} and Gy\"orfi, Morvai, Yakowitz~\cite{GYMY98} it is clear that even for  
the class of all stationary and ergodic  binary 
time-series with almost surely continuous conditional probability 
$P(X_1=1|\dots,X_{-1},X_{0})$ one can not solve Problem~2.

\smallskip
\noindent
For $n\ge 1$, let the function $p_n(\cdot)$   be defined as  
\begin{equation}
\label{defjointdistr}
p_n(x_{-n+1},\dots,x_{0})=P(X_{-n+1}=x_{-n+1},\dots,X_0=x_0)
\end{equation}
where $x_{-i}\in \{0,1\}$ for $0\le i\le n-1$.

\smallskip
\noindent
The entropy rate $H$ associated with a stationary binary time-series
$\{X_n\}^{\infty}_{-\infty}$ is defined as
$
H=\lim_{n\to\infty} -{1\over n} E \log_2 p_n(X_{-n+1},\dots,X_{-1},X_0)$. 
We note that the entropy rate of a stationary binary time-series always exists. 
For details cf. Cover, Thomas \cite{CT91}, pp. 63-64.  

\smallskip
\noindent
Now we may pose our problem.

\smallskip
\noindent
{\bf Problem 3} {\it  
Is there a sequence of strictly increasing stopping times $\{\lambda_n\}$ with 
$$\lambda_n\le 2^{n(H+\epsilon)}$$ 
and an estimation scheme 
$f_n(X_0,\dots,X_{\lambda_n})$ which  depends on the observed data segment 
$(X_0,\dots,X_{\lambda_n})$ such that 
$$
\lim_{n\to\infty} |f_n(X_0,\dots,X_{\lambda_n})-P(X_{\lambda_n+1}=1|X_0,\dots,X_{\lambda_n})|=0
$$
almost surely for all stationary and ergodic binary time series $\{X_n\}_{n=-\infty}^{\infty}$ 
with almost surely continuous conditional probability $P(X_1=1|\dots,X_{-1},X_{0})$?
}

\smallskip
\noindent
It turns out that the  answer  is affirmative and 
such a scheme will be exhibited below.
This result can be interpreted as if the market analyst  can refrain from predicting, that is,
 he may say that he does not want to predict today, but will predict at infinitely many time
 instances, 
and not too rarely, since $\lambda_n\le 2^{n(H+\epsilon)}$,  
and the difference between the prediction and the true conditional probability will 
vanish  almost surely at these stopping times.   
 We note that the stationary processes with almost surely continuous 
conditional distribution  generalize the 
processes for which  the  conditional distribution is actually continuous, 
these are essentially the Random Markov Processes of Kalikow~ 
\cite{Ka90}, or the continuous g-measures studied by Mike Keane in \cite{Ke72}.  
Morvai~\cite{Mo00} proposed a different estimator which is consistent on a certain  stopping time sequence, but those stopping times grow like
an exponential tower which is unrealistic and much faster growth than the mere exponential one in Problem~3.   

\section{The Proposed Estimator}

Let $\{X_n\}_{n=-\infty}^{\infty}$ be a stationary time series taking values from a binary alphabet 
${\cal X}=\{0,1\}$. (Note that all stationary time series $\{X_n\}_{n=0}^{\infty}$ 
can be thought to be a 
two sided time series, that is, $\{X_n\}_{n=-\infty}^{\infty}$. )  
Now we exhibit an estimator which is consistent  on a certain stopping time sequence 
for a restricted class of stationary time series.
For notational convenience, let $X_m^n=(X_m,\dots,X_n)$,
where $m\le n$.

\smallskip
\noindent
Define the stopping times as follows.  Set $\zeta_0=0$. 
For $k=1,2,\ldots$,  define sequence $\eta _k$ and $\zeta_k$
recursively. 
Let
$$
{\eta}_k=
\min\{t>0 : X_{\zeta_{k-1}-(k-1)+t}^{\zeta_{k-1}+t}=X_{\zeta_{k-1}-(k-1)}^{\zeta_{k-1}}\}
\ \ \mbox{and} \ \ 
\zeta_k=\zeta_{k-1}+\eta_k.
$$
One denotes the $k$th estimate of $P(X_{\zeta_k+1}=1|X_0^{\zeta_k})$ by $g_k$, 
and defines it to be
\begin{equation}
\label{fkdistrestimate2}
g_k=
{1\over k}\sum_{j=0}^{k-1} X_{\zeta_j+1}. 
\end{equation} 

\smallskip
\noindent
It will be useful to define other processes 
$\{ {\tilde X}_n\}_{n=-\infty}^{0}$ and 
$\{ {\hat X}^{(k)}_n\}_{n=-\infty}^{\infty}$ for $k\ge 0$ as follows.
Let 
\begin{equation}
\label{defprocesses}
{\tilde X}_{-n}=X_{\zeta_n-n} \ \ \mbox{for $n\ge 0$, and} \ \  
\hat X^{(k)}_{n}=X_{\zeta_k+n} \ \ \mbox{for $-\infty<n<\infty$.}
\end{equation}

\noindent
For an arbitrary stationary binary time series $\{Y_n\}$, and  for all $k\ge 1$ and $1\le i\le k$ 
define $\hat\zeta^k_0(Y^0_{-\infty})=0$ and  
$$
{\hat\eta}^k_i(Y^0_{-\infty})=
\min\{t>0 : 
Y_{\hat\zeta^k_{i-1}-(k-i)-t}^{\hat\zeta^k_{i-1}-t}
=
{Y}_{\hat\zeta^k_{i-1}-(k-i)}^{\hat\zeta^k_{i-1}}\}
$$
and
$$
\hat\zeta^k_i(Y^0_{-\infty})=\hat\zeta^k_{i-1}(Y^0_{-\infty})-\hat\eta^k_i(Y^0_{-\infty}).
$$
When it is obvious on which time series  ${\hat\eta}^k_i(Y^0_{-\infty})$ and  
$\hat\zeta^k_i(Y^0_{-\infty})$
 are evaluated,  
we will use the notation  ${\hat\eta}^k_i$ and  $\hat\zeta^k_i$.
Let $T$ denote the left shift operator,  
that is, $(T x^{\infty}_{-\infty})_i=x_{i+1}$. It is easy to see that if $\zeta_k(x_{-\infty}^{\infty})=l$ then 
${\hat \zeta}^k_k(T^l x_{-\infty}^{\infty})=-l$.

\smallskip
\noindent
We will need the next lemma for later use. 

\begin{lemma} Let $\{X_n\}_{n=-\infty}^{\infty}$ be a stationary binary process. Then 
the  time series $\{{\hat X}^{(k)}_n\}_{n=-\infty}^{\infty}$, 
$\{{\tilde X}_n\}_{n=-\infty}^0$ 
and $\{X_n\}_{n=-\infty}^{\infty}$ have 
identical  distribution. 
Thus all these time series are stationary, and $\{{\tilde X}_n\}_{n=-\infty}^0$ 
can be thought to be two sided stationary time series $\{{\tilde X}_n\}_{n=-\infty}^{\infty}$.  
\end{lemma}
%\begin{pf}

\noindent
Let $k\ge 0$, $n\ge 0$, $m\ge 0$,  
$x^{m}_{m-n}\in {\cal X}^{n+1}$ be arbitrary. 
It is immediate that for $l\ge 0$,  
\begin{equation}\label{shiftequation} 
T^{l} \{X^{\zeta_k+m}_{\zeta_k+m-n}=x^m_{m-n},\zeta_k=l\} =
\{ X^{m}_{m-n}=x^{m}_{m-n},{\hat \zeta}^k_k(X^0_{-\infty})=-l\}. 
\end{equation}

\noindent
First we prove that for $k\ge 0$, $P(({\hat X}^{(k)}_{m-n},\dots,{\hat X}^{(k)}_{m})=(x_{m-n},\dots,x_{m}))=P(X^m_{m-n}=x^m_{m-n})$. 
By the construction in~(\ref{defprocesses}),  
the stationarity of  the time series $\{X_n\}$, and~(\ref{shiftequation}) 
we have 
\begin{eqnarray*}
\lefteqn{ P(({\hat X}^{(k)}_{m-n},\dots,{\hat X}^{(k)}_{m})=(x_{m-n},\dots,x_{m}))}\\
&=&P(X^{\zeta_k+m}_{\zeta_k+m-n}=x^m_{m-n})\\
&=& \sum_{l=0}^{\infty} P(X^{\zeta_k+m}_{\zeta_k+m-n}=x^m_{m-n},\zeta_k=l)\\
&=& \sum_{l=0}^{\infty} P(X^{m}_{m-n}=x^m_{m-n},{\hat \zeta}^k_k(X^0_{-\infty})=-l)\\
&=& P(X^m_{m-n}=x^m_{m-n}).
\end{eqnarray*}

\noindent
Now we prove that  $P({\tilde X}^0_{-n}=x^0_{-n})=P(X^0_{-n}=x^0_{-n})$. 
By the construction in~(\ref{defprocesses}),  
the stationarity of  the time series $\{X_n\}$, and~(\ref{shiftequation}) (with $m=0$)  
we have  
\begin{eqnarray*}
\lefteqn{ P({\tilde X}^0_{-n}=x^0_{-n})=P(X^{\zeta_n}_{\zeta_n-n}=x^0_{-n})}\\
&=& \sum_{l=0}^{\infty} P(X^{\zeta_n}_{\zeta_n-n}=x^0_{-n},\zeta_n=l)\\
&=& \sum_{l=0}^{\infty} P(X^{0}_{-n}=x^0_{-n},{\hat \zeta}^n_n(X^0_{-\infty})=-l)\\
&=& P(X^0_{-n}=x^0_{-n}).
\end{eqnarray*}
The proof of the  Lemma is complete. 
%\end{pf}

\noindent
Now we show the consistency of our estimate $g_k$ defined in~(\ref{fkdistrestimate2}).

\begin{theorem} \label{Theorem3} Let  $\{X_n\}$ be a stationary binary time series. 
For the estimator defined 
in~(\ref{fkdistrestimate2}),
$$
\lim_{k\to\infty} \left| g_k- P(X_{\zeta_k+1}=1|X_0^{\zeta_k})\right| =0\ 
\ \mbox{almost surely}
$$
provided that the conditional probability $P(X_1=1|X_{-\infty}^{0})$ is almost surely 
continuous. Moreover, under the same 
conditions, 
$$
\lim_{k\to\infty} g_k=\lim_{k\to\infty}  
P(X_{\zeta_k+1}=1|X_0^{\zeta_k}) =
P({\tilde X}_1=1|{\tilde X}^{0}_{-\infty}) \ \ \mbox{almost surely.}
$$  
\end{theorem}
%\begin{pf}
Recalling~(\ref{fkdistrestimate2}) we can write
\begin{eqnarray}
\nonumber
g_k &=& {1\over k}\sum_{j=0}^{k-1} [  X_{\zeta_j+1}  -P(X_{\zeta_j+1}=1|X_{-\infty}^{\zeta_j})] + \nonumber {1\over k}\sum_{j=0}^{k-1} P(X_{\zeta_j+1}=1|X_{-\infty}^{\zeta_j})\\
&=& \label{decomposition} {1\over k}\sum_{j=0}^{k-1} \Gamma_j+
{1\over k}\sum_{j=0}^{k-1}  P(X_{\zeta_j+1}=1|X_{-\infty}^{\zeta_j}).
\end{eqnarray}

\noindent
Observe that
 $\{\Gamma_j,\sigma(X_{-\infty}^{\zeta_j+1})\}$ is a bounded martingale difference sequence 
for $0\leq j<\infty$. To see this notice that
$\sigma(X_{-\infty}^{\zeta_j+1})$ 
is monotone
increasing, and  $\Gamma_j$   is measurable with respect to 
$\sigma(X_{-\infty}^{\zeta_j+1})$, and $E(\Gamma_j|X_{-\infty}^{\zeta_{j-1}+1})=0$ for $0\leq j<\infty$ (where you may define $\zeta_{-1}=-1$).
Now apply  Azuma's exponential
bound for bounded martingale differences in Azuma~\cite{Azuma67} to
get that for any $\epsilon>0$,
$$P\left(\left|{1\over k}\sum_{j=0}^{k-1} \Gamma_j\right|>\epsilon\right)
\le 2\exp(-\epsilon^2 k/2).$$
After summing the right hand 
side over $k$, and appealing to
the Borel-Cantelli lemma for a sequence of $\epsilon$'s tending to zero we get  
${1\over k}\sum_{j=0}^{k-1} \Gamma_j \to 0$ almost surely.

\noindent
Define the function $p : {\cal X}^{*-}\rightarrow [0,1]$ as 
$
p(x^{0}_{-\infty})=P(X_1=1|X^{0}_{-\infty}=x^0_{-\infty})
$.

\noindent
For  arbitrary $j\ge 0$, by 
the construction in (\ref{defprocesses}),     
\begin{equation}\label{firstjbitequal}
X_{\zeta_j-j}^{\zeta_j}=
({\hat X}^{(j)}_{-j},\dots,{\hat X}^{(j)}_{0})=
{\tilde X}^0_{-j} \ \ \mbox{and} \ \  
\lim_{j\to\infty} 
d^*({\tilde X}^0_{-\infty},(\dots,{\hat X}^{(j)}_{-1},
{\hat X}^{(j)}_{0}))=0 
\end{equation}
almost surely.
By assumption, the function $p(\cdot)$ is continuous on a set $C\subseteq {\cal
X}^{*-}$ with $P(X^0_{-\infty}\in C)=1$, and  by the Lemma,  
$P({\tilde X}^{0}_{-\infty}\in C)=1$, and  for each $j\ge 0$,
$P((\dots,{\hat X}^{(j)}_{-1},{\hat X}^{(j)}_0)\in C)=1$, and finally,  
$$P({\tilde X}^{0}_{-\infty}\in C, (\dots,{\hat X}^{(j)}_{-1},{\hat X}^{(j)}_0)\in C \ \ \mbox{for all $j\ge 0$})=1.$$  
By the Lemma, 
the construction in~(\ref{defprocesses}),     
the continuity of $p(\cdot)$ on the set $C$, and by~(\ref{firstjbitequal})  
$$
P(X_{\zeta_j+1}=1|X^{\zeta_j}_{-\infty})=p(\dots,{\hat  X}^{(j)}_{-1},{\hat  X}^{(j)}_{0})
\to p({\tilde X}^{0}_{-\infty})=
P({\tilde X}_1=1|{\tilde X}^{0}_{-\infty})
$$
 and 
$
{1\over k}\sum_{j=0}^{k-1}  P(X_{\zeta_j+1}=1|X_{-\infty}^{\zeta_j})\to 
P({\tilde X}_1=1|{\tilde X}^{0}_{-\infty})
$ almost surely. 
We have proved that $g_k\to P({\tilde X}_1=1|{\tilde X}^{0}_{-\infty})$ almost surely.

\noindent
Now observe that by~(\ref{defdistance}) and 
the continuity of $p(\cdot)$ on the set $C$, almost surely, for all $\epsilon>0$, there is a   
$J(\epsilon,{\tilde X}^0_{-\infty})$, such that for all $z^0_{-\infty}\in C$, if $z^0_{-J}={\tilde X}^0_{-J}$ then 
$
|p(z^0_{-\infty})-p({\tilde X}^0_{-\infty})|<\epsilon
$.
By~(\ref{firstjbitequal}), 
and since $\epsilon>0$ was arbitrary, almost surely, 
\begin{eqnarray*}
\lim_{j\to\infty} 
P(X_{\zeta_j+1}=1|X^{\zeta_j}_{0})&=&
\lim_{j\to\infty} 
E\{ P(X_{\zeta_j+1}=1|X^{\zeta_j}_{-\infty}) |X^{\zeta_j}_0\}\\
&=&
\lim_{j\to\infty} 
E\{ p(X^{\zeta_j}_{-\infty}) |X^{\zeta_j}_0\}\\
&=& p({\tilde X}^0_{-\infty})=P({\tilde X}_1=1|{\tilde X}^{0}_{-\infty}).
\end{eqnarray*}
The proof of Theorem~\ref{Theorem3} is complete.
%\end{pf}

\noindent
{\bf Remark.} We note that for all stationary binary time-series, the estimation scheme described above is consistent in probability. This may
be seen as follows:    
\begin{eqnarray*}
\lefteqn{  E\left| g_k-
P(X_{\zeta_k+1}=1|X_{0}^{\zeta_k})\right| }\\
&\le&
 E\left| {1\over k}\sum_{j=0}^{k-1} [  X_{\zeta_j+1}  
-
P(X_{\zeta_j+1}=1|X_{-\infty}^{\zeta_j})] \right| \\
&+&  {1\over k}\sum_{j=0}^{k-1}E\left| 
P({\hat X}^{(j)}_{1}=1|\dots,{\hat X}^{(j)}_{-1},{\hat X}^{(j)}_0)-
P({\hat X}^{(j)}_{1}=1|{\hat X}^{(j)}_{-j},\dots,{\hat X}^{(j)}_0)
\right|\\
&+&
  E\left| {1\over k}\sum_{j=0}^{k-1}
P({\hat X}^{(k)}_{1}=1|{\hat X}^{(k)}_{-j},\dots,{\hat X}^{(k)}_0)-
P({\hat X}^{(k)}_{1}=1|{\hat X}^{(k)}_{{\hat\zeta}^k_k},\dots,{\hat X}^{(k)}_0)
\right|,\\
\end{eqnarray*}
where we used~(\ref{firstjbitequal}) and the Lemma.
The first term converges to zero since 
$X_{\zeta_j+1}-P(X_{\zeta_j+1}=1|X_{-\infty}^{\zeta_j})$ 
is a martingale difference sequence with respect to 
$\sigma(X_{-\infty}^{\zeta_j+1})$ 
and an average of bounded martingale differences 
converges to zero almost surely cf. Azuma~\cite{Azuma67}. 
Applying~(\ref{defprocesses}), (\ref{firstjbitequal}) and the Lemma, the sum of the  last two terms can be  estimated by the sum     
\begin{eqnarray*}
\lefteqn{
 {1\over k}\sum_{j=0}^{k-1} E\left|P(X_1=1|X_{-\infty}^0)-
P(X_1=1|X_{-j}^0)
\right|}  \\
&+&  E\left|{1\over k}\sum_{j=0}^{k-1} P(X_1=1|X_{-j}^0)-P(X_1=1|X_{{\hat\zeta}^k_k}^0)\right|
\end{eqnarray*}
and both terms  converge to zero since by the martingale convergence theorem
$
\lim_{j\to\infty} P(X_1=1|X_{-j}^0)=P(X_1=1|X_{-\infty}^0)
$ almost surely, and thus the limit in fact exists and equals zero.

\smallskip
\noindent
Next we will give some universal estimates for the growth rate of the stopping times $\zeta_k$ in terms of the entropy rate of the process. This
is natural since the $\zeta_k$ are defined by recurrence times for blocks of length $k$, and these are known to grow exponentially with the
entropy rate. (Cf. Ornstein and Weiss~\cite{OrWe93}.)  

\begin{theorem} \label{Theorem4}
Let  $\{X_n\}$ be a stationary and ergodic binary time series.  
Then for arbitrary $\epsilon>0$, 
$$
\zeta_k< 2^{k(H+\epsilon)} 
\ \ \mbox{ eventually almost surely,} 
$$
where $H$ denotes the entropy rate associated with time series $\{X_n\}$. 
\end{theorem}
%\begin{pf}

\noindent 
Let ${\cal X}^*$ be the set of all two-sided  binary sequences, that is, 
$${\cal X}^*= \{ (\dots,x_{-1},x_0,x_1,\dots):
x_i\in \{0,1\} \ \ \mbox{for all $-\infty<i< \infty$}\}.
$$

\noindent
Define $B_k\subseteq \{0,1\}^k$ as 
$$
B_k=\{ x^0_{-k+1} \in \{0,1\}^k: 2^{-k(H+0.5\epsilon)} < p_k(x^0_{-k+1})\}, 
$$
where $p_k(\cdot)$ is as in~(\ref{defjointdistr}).
Note that there is a trivial bound on the cardinality of the set $B_k$, namely,
\begin{equation} \label{cardbound}
|B_k|\le 2^{k(H+0.5\epsilon)}. 
\end{equation}
By the Lemma, the distribution of the time series 
$\{\tilde X_n\}$ is the same as the distribution of $\{X_n\}$ 
and  by the Shannon-McMillan-Breiman Theorem (cf. Cover, Thomas \cite{CT91}, p. 475),  
\begin{equation} \label{SMcMillan}
P\left( \bigcup_{k=1}^{\infty}\bigcap_{i\ge k} \{ {\tilde X}^0_{-i+1}\in B_i\}\right) =1.
\end{equation}

\noindent
Define the set $Q_k(y^0_{-k+1})$ as follows: 
$$
Q_k(y^0_{-k+1})=\{z^{\infty}_{-\infty}\in {\cal X}^* :  -{\hat \zeta}^k_k(z^0_{-\infty}) 
\ge 2^{k(H+\epsilon)}, z^0_{-k+1}=y^0_{-k+1})\}.
$$
We will estimate the probability of $Q_k(y^0_{-k+1})$  by means of the ergodic theorem. 
Let $x^{\infty}_{-\infty}\in {\cal X}^*$ be a typical sequence of the time series $\{X_n\}$. 
Define $\alpha_0(y^0_{-k+1})=0$ 
and for $i\ge 1$ let 
$$
\alpha_i(y^0_{-k+1})=\min \{l> \alpha_{i-1}(y^0_{-k+1}): T^{-l} x_{-\infty}^{\infty}\in Q_k(y^0_{-k+1})\}.
$$
Define also $\beta_0(y^0_{-k+1})=0$  
and for $i\ge 1$ let 
$$
\beta_i(y^0_{-k+1})=\min \{l> \beta_{i-1}(y^0_{-k+1})+2^{k(H+\epsilon)}: 
T^{-l} x_{-\infty}^{\infty}\in Q_k(y^0_{-k+1})\}.
$$
Observe that for arbitrary $l>0$, 
$$
\sum_{j=1}^{\infty}   1_{\{\beta_{l-1}(y^0_{-k+1}) < 
\alpha_j(y^0_{-k+1}) \le \beta_l(y^0_{-k+1})\}}\le k+1.
$$
By the Lemma and the ergodicity of the time series $\{X_n\}$, 
\begin{eqnarray}
\lefteqn{ \nonumber
P((\dots, {\hat X}^{(k)}_{-1}, {\hat X}^{(k)}_0, {\hat X}^{(k)}_{1},\dots )\in Q_k(y^0_{-k+1}))=
P( X_{-\infty}^{\infty} \in Q_k(y^0_{-k+1}) ) }\\
&=& \nonumber
\lim_{t\to\infty} {1\over \beta_t(y^0_{-k+1})} \sum_{j=1}^{\infty} 
1_{ \{\alpha_j(y^0_{-k+1})\le \beta_t(y^0_{-k+1})\} }  \\
&=&  \nonumber\lim_{t\to\infty}{1\over \beta_t(y^0_{-k+1})} 
\sum_{l=1}^{t} \sum_{j=1}^{\infty}   1_{\{\beta_{l-1}(y^0_{-k+1}) < 
\alpha_j(y^0_{-k+1}) \le \beta_l(y^0_{-k+1})\}} \\
&\le&  \lim_{t\to\infty} {t (k+1)\over t 2^{k(H+\epsilon)} }
= \label{upperbound} {(k+1)\over 2^{k(H+\epsilon)}}. 
\end{eqnarray}
By the construction in (\ref{defprocesses}),   
$-{\hat \zeta}_k^k(\dots,{\hat X}^{(k)}_{-1},{\hat X}^{(k)}_{0})=\zeta_k(X_0^{\infty})$, and 
$({\hat X}^{(k)}_{-k+1},\dots,{\hat X}^{(k)}_0)={\tilde X}^0_{-k+1} $ 
and by the upper bound on the cardinality of set $B_k$   
in (\ref{cardbound}) and by (\ref{upperbound}), 
   we get  
\begin{eqnarray*}
\lefteqn{
P(\zeta_k(X_0^{\infty})\ge 2^{k(H+\epsilon)}, {\tilde X}^0_{-k+1}\in B_k)}\\
&=&
P(-{\hat \zeta}_k^k(\dots,{\hat X}^{(k)}_{-1}, {\hat X}^{(k)}_{0}) 
\ge 2^{k(H+\epsilon)}, {\tilde X}^0_{-k+1}\in B_k)\\
&=& 
P(-{\hat \zeta}_k^k(\dots,{\hat X}^{(k)}_{-1}, {\hat X}^{(k)}_{0}) 
\ge 2^{k(H+\epsilon)}, 
({\hat X}^{(k)}_{-k+1},\dots,{\hat X}^{(k)}_0) \in B_k)\\
&=& \sum_{y^0_{-k+1}\in B_k} 
P((\dots, {\hat X}^{(k)}_{-1}, {\hat X}^{(k)}_0, {\hat X}^{(k)}_{1},\dots )\in Q_k(y^0_{-k+1}))
\le 
(k+1) 2^{-k 0.5\epsilon}.
\end{eqnarray*}
The right hand side sums, the Borel-Cantelli Lemma and the Shannon-McMillan-Breiman Theorem 
in (\ref{SMcMillan}) together yield that  
$\zeta_k< 2^{k(H+\epsilon)}$eventually almost surely
and Theorem~\ref{Theorem4} is proved. 

%\end{pf}

% -------  T E X T  E N D S  -----------------------------

\end{document}